\documentclass[graybox]{svmult}

\usepackage{url}

\usepackage{type1cm}        
\usepackage{makeidx}         
\usepackage{graphicx}        
\usepackage{multicol}        
\usepackage[bottom]{footmisc}

\usepackage{newtxtext}       %
\usepackage{newtxmath}       

\makeindex             

\begin{document}

\title*{A data analysis of women's trails among ICM speakers}
\author{Helena Mihaljevi\'{c} and Marie-Fran\c{c}oise Roy}
\institute{Helena Mihaljevi\'{c} \at Hochschule f\"{u}r Technik und Wiertschaft Berlin,  Wilhelminenhofstra{\ss}e 75A,
12459 Berlin, Germany,  \email{Helena.Mihaljevic@HTW-Berlin.de}
\and Marie-Fran\c{c}oise Roy \at IRMAR (UMR CNRS 6625), Universit\'{e} de Rennes 1, Campus de Beaulieu, 35042 Rennes Cedex, France, \email{marie-francoise.roy@univ-rennes1.fr}}
%
%
\maketitle

\abstract*{The International Congress of Mathematicians (ICM), inaugurated in 1897, is the greatest effort of the mathematical community to strengthen international communication and connections across all mathematical fields. Meetings of the ICM have historically hosted some of the most prominent mathematicians of their time. Receiving an invitation to present a talk at an ICM signals the high international reputation of the recipient, and is akin to entering a `hall of fame for mathematics'. Women mathematicians attended the ICMs from the start. With the invitation of Laura Pisati to present a lecture in 1908 in Rome and the plenary talk of Emmy Noether in 1932 in Zurich, they entered the grand international stage of their field. At the congress in 2014 in Seoul, Maryam Mirzakhani became the first woman to be awarded the Fields Medal, the most prestigious award in mathematics. In this article, we dive into assorted data sources to follow the footprints of women among the ICM invited speakers, analyzing their demographics and topic distributions, and providing glimpses into their diverse biographies.}

\abstract{The International Congress of Mathematicians (ICM), inaugurated in 1897, is the greatest effort of the mathematical community to strengthen international communication and connections across all mathematical fields. Meetings of the ICM have historically hosted some of the most prominent mathematicians of their time. Receiving an invitation to present a talk at an ICM signals the high international reputation of the recipient, and is akin to entering a `hall of fame for mathematics'. Women mathematicians attended the ICMs from the start. With the invitation of Laura Pisati to present a lecture in 1908 in Rome and the plenary talk of Emmy Noether in 1932 in Zurich, they entered the grand international stage of their field. At the congress in 2014 in Seoul, Maryam Mirzakhani became the first woman to be awarded the Fields Medal, the most prestigious award in mathematics. In this article, we dive into assorted data sources to follow the footprints of women among the ICM invited speakers, analyzing their demographics and topic distributions, and providing glimpses into their diverse biographies.}

\section{The hall of fame for mathematics}
\label{sec:1}
Ever since its inaugural gathering in 1897, the International Congress of Mathematicians (ICM) has signified the greatest effort of the mathematical community to establish international communication and connection across all mathematical topics. Throughout their history, the congresses have hosted some of the most prominent mathematicians of their time. Needless to say, receiving an invitation to present a talk at an ICM is a matter of high international reputation, often compared with the entrance into a `hall of fame for mathematics'. In fact, it is no exaggeration to state that an ICM invitation is often treated like the reception of a major research award.

Women mathematicians attended the ICM from the start, not only as accompanying persons but also participating on their own, e.g. as professional mathematicians.  Nevertheless, female speakers remained very few. The share of women in selected congresses has been addressed in some previous works. Fulvia Furinghetti studied the presence and contribution of women to the discipline of mathematics education in the first half of the 20th century using data from two scientific journals, the proceedings of the first International Congress on Mathematical Education (ICME) in 1969 and the didactics sections of the ICM proceedings until 1966 \cite{Fur}. She describes the difficulties posed by differences in structure of the individual congresses, the layout of the congress proceedings, or inconclusive and incomplete data (e.g. regarding the distinction of `accompanying persons'). She provides numbers of women among participants and contributors for ICMs until 1966 and gives insights into biographies of women pioneers. The essays \cite{Sad} by Cora Sadosky and \cite{CaL} by Bettye Anne Case and Anne M. Leggett from the collection `Complexities. Women in Mathematics' address the participation of women lecturers since 1974, focussing mainly on the collective efforts of women in the 1970s and 1980s to overcome their persistent underrepresentation as invited congress speakers. Both pieces arrive at similar conclusions, namely that the actions in the 1970s and 1980s have strongly contributed to the diversification of the congress, yielding a significantly higher chance for qualified women to be invited to speak.

While the mentioned research addresses the participation of women at individual congresses or throughout certain periods, to our knowledge there is no global exploratory analysis of the demographics of ICM speakers from its beginning until today. In this contribution, we thus investigate data on all invited ICM speakers from 1897 to 2018. Using various data sources, in particular the list of all invited speakers from Wikipedia, Wikidata pages of individual speakers, and the subdivision of congress speakers into sections from the International Mathematical Union (IMU), we are able to address the following questions regarding women's participation: How inclusive has the congress been throughout its history? What factors might have positively influenced the share of women? Are there noteworthy differences between women and men speakers regarding age, country of residence or research areas?

We start out by describing pioneer contributions of women. We then outline the development of women's participation in the congresses over time, elaborating some of the advances and setbacks. Finally, we investigate the distribution of women and men speakers by countries of citizenship and sections of the delivered talk.

\section{Data basis and methods}
\label{sec:2}
In February 2018, we programmatically extracted all names and ICM dates from the Wikipedia website `List of International Congresses of Mathematicians Plenary and Invited Speakers' \cite{Wik}, which resulted in a table containing 3,745 plenary and invited speakers from 1897 until 2014. Using the hyperlinks contained therein, we retrieved the gender, country of citizenship, date of birth and employer from Wikidata, a free, human- and machine-readable knowledge database that serves as a central storage for structured data of other Wikimedia projects, including Wikipedia\footnote{Every Wikipedia article is supposed to have a corresponding entry in Wikidata.}. We have found a Wikidata page for 82.6\% of all listed speakers, and 77.5\% of all unique individuals (various mathematicians gave multiple talks). A Wikidata page existed for almost all women, namely 92.5\%. The coverage shows a certain trend: with exception of the large congresses in 1928 and 1932, the vast majority of speakers from the early congresses (usually 90\% or higher) has a page in Wikidata, with decreasing trend over time. 

For speakers invited to the ICM 2018 in Rio de Janeiro, we extracted their names, country of citizenship and the ICM sections of their talks from the official ICM-2018 website. 
We used Python package gender-guesser\footnote{https://pypi.org/project/gender-guesser/}, which has shown very reliable results in a recent benchmark on name-based gender inference \cite{SaM}, to infer the gender\footnote{For all authors we used a gender assignment provided by a third party (Wikidata or a web service) which, for our dataset, resulted in a binary schema.} of the speakers using their forenames when this information was missing. For speakers whose names are not highly correlated with only one gender (across different countries and languages), and for which gender-guesser hence did not produce a definite gender assignment, we filled this information manually, mainly based on field knowledge and Internet research.

The International Mathematical Union (IMU) provided us with a file containing speaker names, ICM date and place, and the name of the section of the corresponding talk (see \cite{ICM} for a search interface within the official IMU website). Due to different name spellings in the datasets, we applied fuzzy string matching techniques to combine the data sources and add the sections to our original data set. For many speakers at the congress in 1950 the section was missing and hence needed to be filled manually.

In addition, we added the date of birth and country of citizenship for all women speakers in order to create a data basis which is as complete as possible and that can be used for information and teaching purposes beyond this analysis. For this purpose we have contacted those women in our list for whom information was still missing. We have made the list of all speakers available at \cite{Mih}.
We noted that there exist different countings of invited ICM speakers. For instance, the list of speakers provided by the IMU \cite{ICM} contains around 400 speakers more than the list at Wikipedia \cite{Wik}, in particular for the congresses before 1950. This is mainly due to the change of terminology over time and the respective counting schemes. On the other hand, the list at Wikipedia contains speakers who were invited but did not attend. Our analyses are based on the Wikipedia list \cite{Wik} which applies the post World War II terminology in which the one-hour speakers in the morning sessions are called `Plenary Speakers' and the usually more numerous speakers (in the afternoon sessions) whose talks are included in the ICM published proceedings are called `Invited Speakers''' \cite{Wik}. Usually, there were a lot more additional shorter contributions that were not always part of the congress  proceedings. Moreover, the list of speakers from Wikipedia \cite{Wik} does not reflect whether a speaker gave more than one talk at a given congress. This, in fact, was not so rare; for example at the congress in 1900 in Paris, G\"{o}sta Mittag-Leffler gave both a plenary talk and one in the \textit{Analysis} section. In order to  take into account such multiple contributions, we have expanded the data using the sections from the list supplied by the IMU \cite{ICM}.

\section{Women pioneers}

The organizers of the congress in 1908 in Rome invited Laura Pisati,  the first woman to present a paper. Not much is known about her personal and professional life, other than that Pisati was an active mathematics researcher and the author of internationally recognized publications. In zbmath.org, we find a book and three research articles listed in her author profile \cite{Pis}, two of them published in the influential \textit{Rendiconti del Circolo Matematico di Palermo}, the journal of the Mathematical Circle of Palermo, of which she was a member. In 1905, she also became a member of the German Mathematical Society \cite[p.12]{Jar}. According to her membership information, Pisati was born in Ancona (date of birth not listed). She graduated in mathematics from the University of Rome in 1905 \cite{Sci}. Since 1897 she had worked as a teacher at the Technical School `Marianna Dionigi' in Roma (Scuola Tecnica `Marianna Dionigi' di Roma), one of the the first secondary schools for girls in Rome. She was engaged to Giovanni Giorgi, an Italian physicist and electrical engineer. In 1900, Pisati had been entrusted with the supervision of his thesis in Mathematics \cite{GGi}. Sadly, she died young on March 30 1908, only a few days before the 1908 congress in Rome and before her planned wedding to Giorgi. Her paper `Saggio di una teoria sintetica delle funzioni di variabile complessa' was presented by a male colleague. 

In the report on the sectional meetings of the congress \cite{Moo}, Laura Pisati appears as the only speaker with first and last name listed, showing the singularity of women's presence in this circle at that time. Interestingly, Giorgi himself was an invited speaker at three subsequent ICMs, in 1924 in Toronto, in 1928 in Bologna, and in 1932 in Zurich. He cited Pisati's work in his 1924 ICM contribution with the words ``See also some very striking results given by LAURA PISATI in her paper Sulle operazioni funzionali non analitiche originate da integrali definiti. Rend. Cire. Mat. Palermo, Tomo XXV (1908) pp. 272-282.'' \cite[p.45]{Fie}. 

Four years later, in 1912, Hilda Hudson was the first woman to speak at an ICM with a paper she presented in the \textit{Geometry} section. Hudson, a member of a family of distinguished mathematicians, worked mainly in the theory of Cremona transformations, on which she had published various articles. Between 1910 and 1913, she was an Associate Research Fellow at the Newnham College\footnote{Newnham College, founded in 1871, was the second women's college to be established in Cambridge.  It acquired full university status in 1948, the year in which the first women were officially admitted to the University.} \cite{BGG}. As pointed out in \cite{Fur}, Hilda Hudson is listed in the Proceedings of the congress in 1912 as an accompanying person to her father, Prof. William Henry Hoar Hudson, showing how misleading the distinction between accompanying persons and 'real' participants  was in that period.

\begin{figure}[b]
\includegraphics[width=\textwidth]{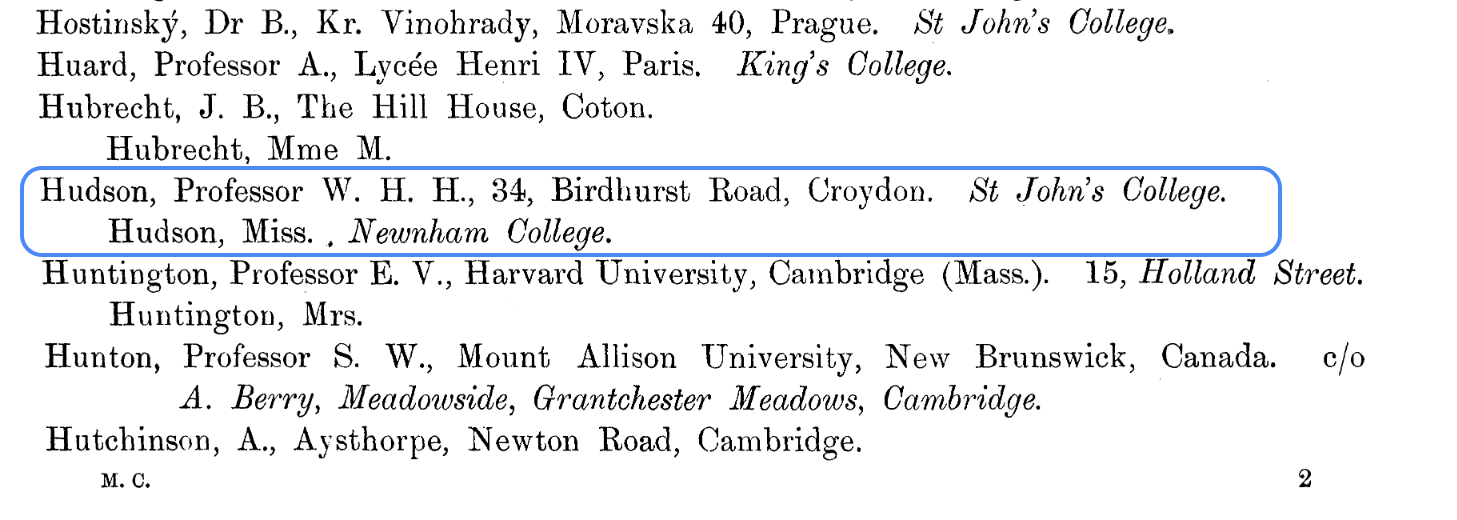}
\caption{Hilda Phoebe Hudson, the first woman who presented her work at an ICM, listed as an accompanying participant at the ICM 1912 in Cambridge.}
\label{fig:hudson}      
\end{figure}

The 1932 congress in Zurich witnessed the first plenary talk by a woman, given by Emmy Noether, who spoke about hypercomplex systems in their relations with commutative algebra and number theory\footnote{Original title of the talk in German: "Hyperkomplexe Systeme in ihren Beziehungen zur kommutativen Algebra und zur Zahlentheorie".}. Her invitation certainly marked a milestone in the representation of women within the international mathematical community. Noether had already attended previous congresses. At the age of 26 she accompanied her father, Max Noether, who spoke at the congress in 1908 in Rome, where Pisati was supposed to present her work. Prior to her plenary lecture in 1932, Emmy Noether gave a talk at the congress in Bologna four years earlier. 
As the positive trend in the early years of the ICM did not persist, it was almost 60 years until Karen Uhlenbeck became the second woman to give a plenary talk at an ICM under the title `Applications of non-linear analysis in topology'.

\begin{figure}[b]
\includegraphics[height=.25\textheight]{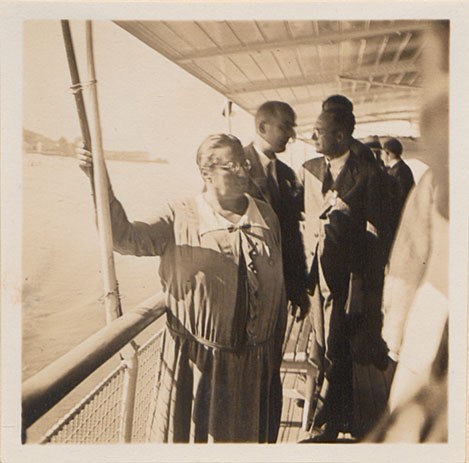}\hfill\includegraphics[height=.25\textheight]{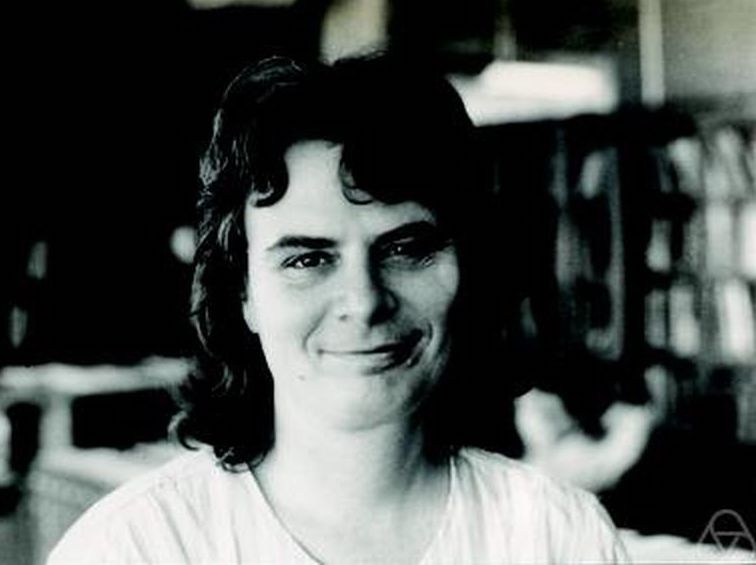}
\caption{Left: Emmy Noether (front) on a steamboat trip during ICM-1932. (ETH-Library, Zurich). Right: Karen Uhlenbeck in 1982, eight years before her plenary lecture in Kyoto. (Oberwolfach Photo Collection)}
\label{fig:noether_uhlenbeck}      
\end{figure}

In 2014 at the ICM in Seoul, Maryam Mirzakhani was awarded the Fields Medal for ``her outstanding contributions to the dynamics and geometry of Riemann surfaces and their moduli spaces''\footnote{ICM laudation, http://www.icm2014.org/en/awards/prizes/f4}.  She is the only woman among the 60 mathematicians who have received
 the Fields Medal, a prize conferred since 1936 to at most four mathematicians at each congress under the age of 40. Mirzakhani was diagnosed with breast cancer in 2013 and died on July 14, 2017, at the age of 40. 

\section{The history in numbers: advances and setbacks}

\begin{figure}[b]
\centering
\includegraphics[width=\textwidth]{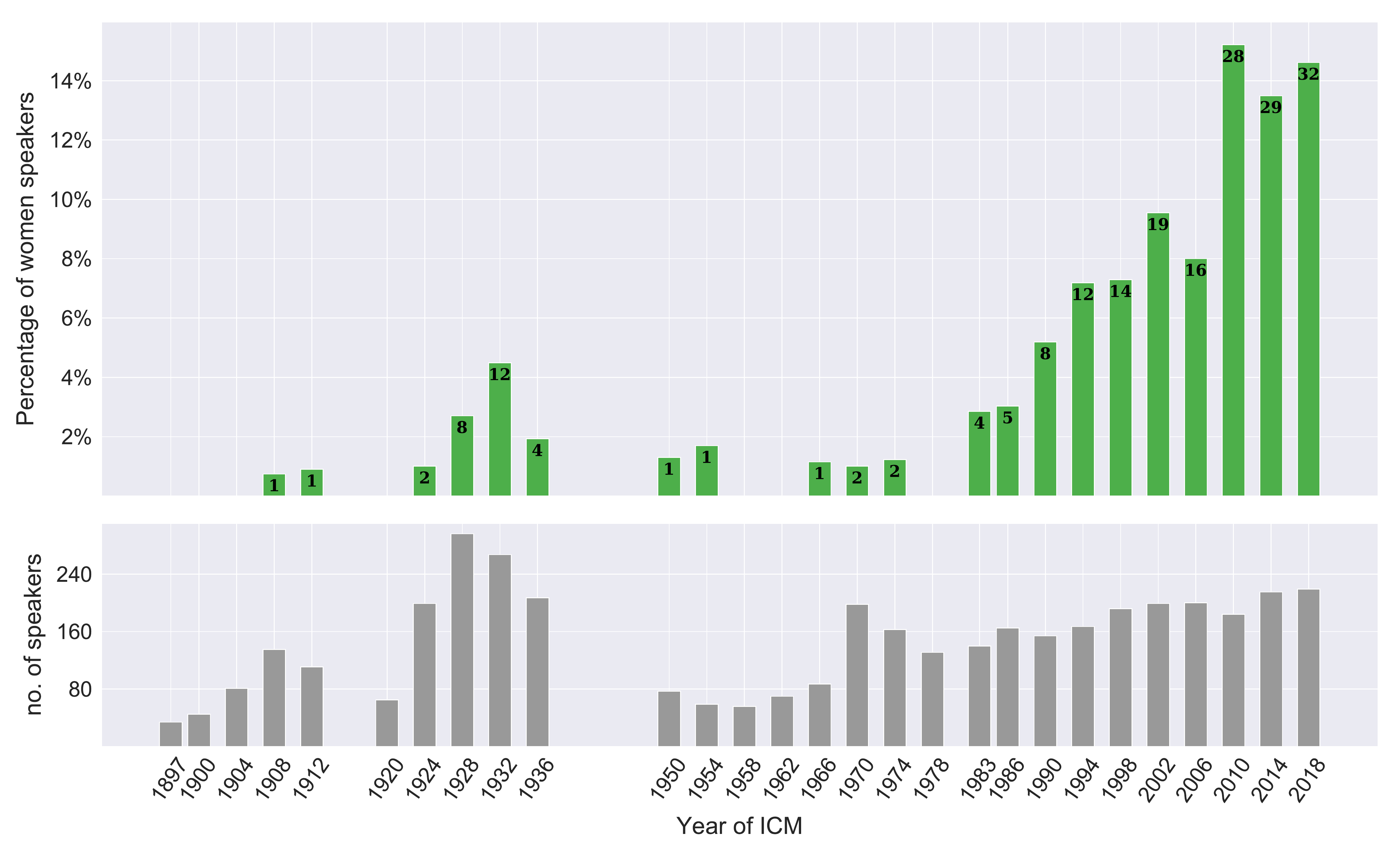}
\caption{Upper panel: Bar height shows the percentage of women speakers per ICM, the numbers inside correspond to the total numbers of invited women. Lower panel: total number of speakers per year.}
\label{fig:longitudinal_dev}      
\end{figure}

Out of 4,120 invited contributions from 1897 to 2018, 202 were presented or authored by women, which amounts to only 5\% of the total. Women's participation over time, however, did not grow steadily but, instead, shows multiple trends. As presented in Fig. \ref{fig:longitudinal_dev}, a comparatively large number of women presented their research at the congresses in 1928 at Bologna and in 1932 at Zurich. This reflects the overall progressive societal and political spirit of the 1920s, which had also enhanced the situation of women in science. In fact, the ICM in 1932 marks a pinnacle in the history of ICMs regarding the role of women. Emmy Noether gave the first plenary lecture by a woman; various women's colleges and organizations of university women sent delegates, among them the Bedford College for Women (London), Hunter College (New York), the International Federation of University Women, and the American Association of University Women. 

ICMs were always affected by global political events. The first substantial tension occurred in the aftermath of World War I, as mathematicians from Germany were excluded during the ICMs in 1920 and 1924. The sole choice of Strasbourg as the location for the congress in 1920 was a political statement in itself. Already in the 1920s, Italy, which experienced a golden era in both pure and applied mathematics at the turn of the 19th century, showed the first signs of a deep crisis caused by a spreading fascism\cite{Isr}. The last congress before World War II that took place in 1936 in Oslo was signified by different political agendas, in particular by the German strategy to present `Aryan mathematics'. Italian mathematicians boycotted the congress, Soviet mathematicians were denied the travel permission by their political authorities. That only few women were invited to the congress in Oslo seems not too surprising given the political situation at that time -- the spread of Fascism through Europe, persecution of Jewish mathematicians, and the worldwide economic depression -- which in some aspects affected women mathematicians on a larger scale. When the first congress after World War II took place in 1950 in Cambridge (USA), the only woman lecturer was Mary Cartwright, at that time Mistress of the Girton College\footnote{Girton College was the first women's college to be established in Cambridge.  It began in Hitchen (about 24 miles from Cambridge) in 1869 before moving to Girton in 1873, when it acquired the name Girton College. Like Newnham, it obtained full college status in 1948.}.

It took 60 years to reach a share of women among ICM speakers comparable to that in 1932. Among the manifold reasons for this situation are undeniably the impact of some historical and political developments. The aftermath of World War II was characterized by a rollback in society as a whole. The 1950s experienced a return to conservative gender roles, in which women were expected to take care of the domestic sphere, leaving the work places to the men who were coming back from the battlefields. These conceptions had impact on university education as well. During the conservative post-war era in Germany, for instance, the share of female students decreased significantly, and there was general agreement that men should take precedence in accessing the limited study places. However, some countries managed to overcome some of these barriers in women's university education and research faster than others. Partially, these general trends are also reflected in country-based differences regarding the presence of women speakers at postwar  ICMs: in the 11 congresses between 1950 and 1990, of the 24 talks given by women, almost all delivered by speakers from the United States, France, United Kingdom, or Russia but none by speakers from Italy or Germany. By contrast, in the ten congresses before World War II of a comparable total of 27 talks by women, three of those speakers were from Germany and four from Italy. 

The situation for women as active participants in ICMs changed in the 1990s and has shown a certain level of stability ever since. In particular, the recent three congresses have witnessed a hitherto unseen participation of women: of all lectures delivered by women in the history of the congress, 80\% took place since the meeting in 1990 in Tokyo. The drastic change affects plenary sessions in particular: The ICMs in 2002 in Beijing, 2010 in Hyderabad and 2018 in Rio de Janeiro collectively accounted for ten of the total 18 plenary lectures by women since the premiere by Emmy Noether in 1932.

Despite the overall progress towards gender equality in mathematics in the recent decades, the increase of women speakers since 1990 cannot be interpreted simply as a positive side effect of a global development. A closer look at the events during the congresses shows that the increased invitation of women speakers is also, and maybe above all, the result of interventions by groups and individuals at various levels. As described in \cite{Sad,CaL}, since 1974, organizations of women such as the Association for Women in Mathematics (AWM) have set up events during the congresses, often sparking discussions on what was often perceived as a systematic omission of women as invited speakers. At various congresses in the 1970s and 1980s  resolutions were passed with the aim to increase the number of lectures by women. At ICM-1974, concerns about the small number of women speakers were raised during a discussion by the AWM. At the next ICM in 1978, a public protest initiated by AWM members resulted in a widely supported resolution to improve the situation of women in the future. Four women were invited to the congress in Warsaw in 1983, but there were no protests or reminders to keep improving the situation. It is probably no coincidence that in the program announcement of ICM-1986, not a single woman was listed in traditional mathematics research areas, suggesting that, as formulated by Sadosky in \cite{Sad}, ``when there are no reminders about women mathematicians, colleagues tend not to remember us''. The program of the ICM-1986 was changed on short notice, again through intervention, by presenting 25 qualified women candidates to the Executive Committee. The informal panel discussion organized by AWM on the situation of women in mathematics that took place during the congress in 1986 was at the origin of the constitution of the European Women in Mathematics (EWM).

The engagement of Mary Ellen Rudin in her role as the head of the U.S. delegates at the IMU General Assembly in 1986 is an illustrative example of what can change when individuals in prominent positions pursue this topic. The president of the ICM-1990 in Kyoto explicitly stated that the committees have followed Rudin's recommendation that subfields of mathematics, women, and mathematicians in small countries should not be overlooked \cite{Sat}. 

Since 2010, specific satellite meetings of the congresses have been organized with the goal of highlighting the contributions and achievements, but also to address concerns of women mathematicians: the International Conference of Women Mathematicians in 2010 in Hyderabad and the International Congress for  Women Mathematicians in 2014 in Seoul. In Rio de Janeiro, the World Meeting for Women in Mathematics (WM\^{}2) was set up as a satellite meeting combined with a panel discussion on the gender gap in the mathematical and natural sciences, and integrated into the ICM-2018 program and the ICM proceedings.

\section{Any difference?}

Within the group of ICM speakers, men and women do show some differences regarding certain demographic aspects. For instance, both women and men speakers were around 44 years old when invited to give a lecture. However, before ICM-1950, women speakers were on average 36 years old, 9 years younger than their male colleagues. Since 1950, their average age has surpassed men's by almost 5 years. 

We have focused on two particular aspects: the country of citizenship and the sections in which the speakers presented their research. The country of citizenship is interesting demographic information for ICMs, in particular due to their regional focus. The mathematical research fields, on the other hand, are known to show high variance in the share of women \cite{MST}.

\subsection{Distribution by countries}
We have collected the country of citizenship for 3,038 out of 3,987 speakers, mostly through their Wikidata pages. Further, we have undertaken additional manual efforts in order to collect missing information for all women speakers by using their websites, personal contacts or contacting them by e-mail. 

Various countries listed in Wikidata do not exist anymore. We have thus replaced their names with those of today's states, e.g. Second Polish Republic with its successor state Poland or the Weimar Republic with Germany. When such a replacement is not possible, in particular for states that have disintegrated over the course of time such as Kingdom of Yugoslavia,  Czechoslovak Socialist Republic or Austria-Hungary, we have inspected the demographics of the corresponding speakers to assign the closest state existing today. 

Furthermore, for speakers listed with more than one country of citizenship we have weighted each of them by inverse frequency, e.g. for a speaker with citizenships of Germany and the United States, each would be counted as one half\footnote{The Wikidata entry of some speakers shows quite a few different citizenships, e.g. \DJ uro Kurepa, a plenary speaker in 1954 and 1958, has had 5 different citizenships according to his Wikidata entry: Socialist Federal Republic of Yugoslavia; Kingdom of Yugoslavia; Kingdom of Serbs, Croatians and Slovenes; Austria-Hungary; Federal Republic of Yugoslavia.}. This procedure was necessary since we do not have information on the time period in which a citizenship was valid. This aggregation thus presumably shows a certain bias for countries which are known to have attracted mathematicians, in particular the United States \cite{And}.

\begin{figure}[b]
\centering
\includegraphics[width=\textwidth]{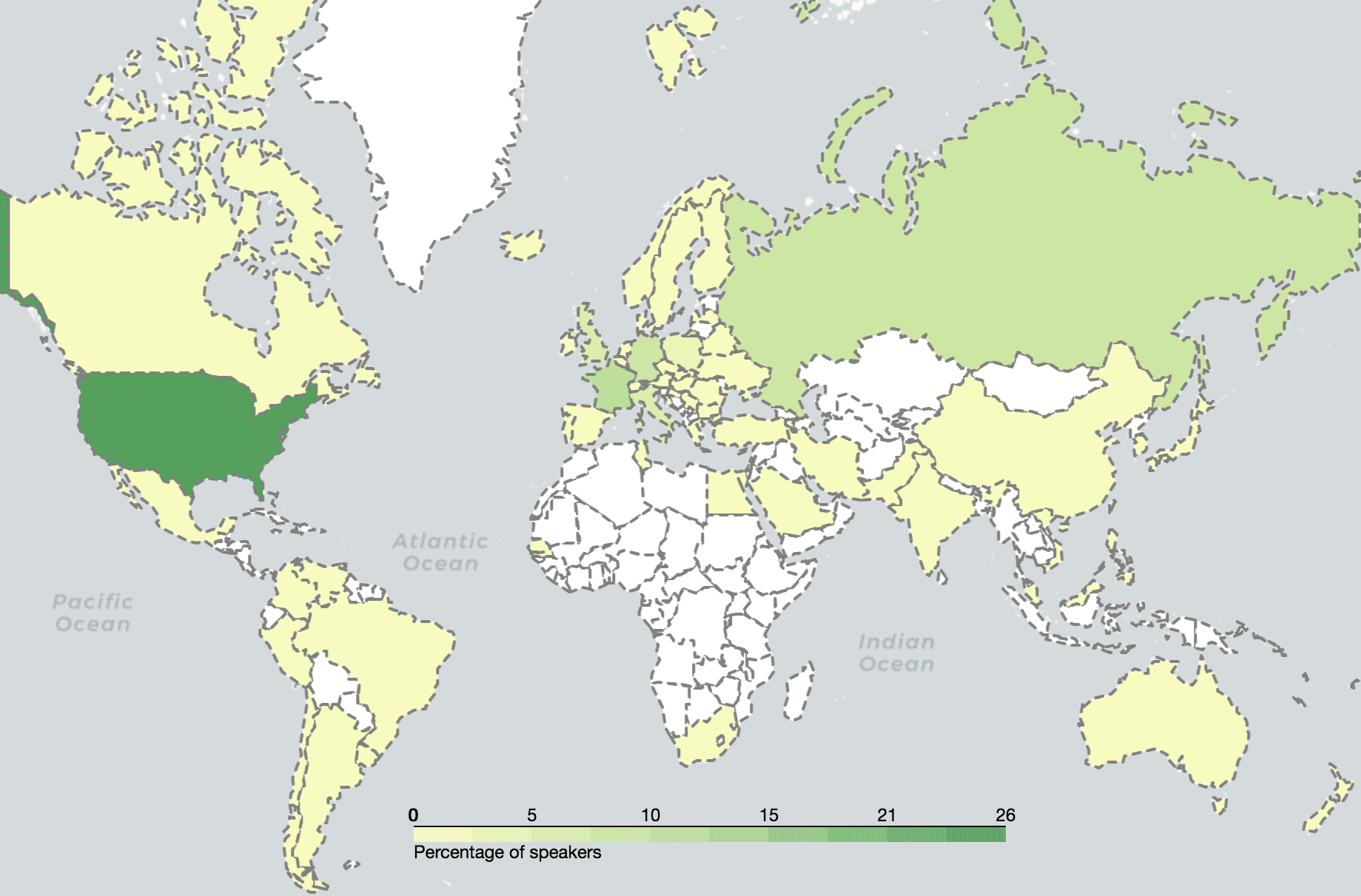}
\caption{Geographical distribution of all women speakers according to their country of residence. }
\label{fig:world_map}      
\end{figure}

The map in Fig. \ref{fig:world_map} shows the proportion of countries as countries of citizenship among all women speakers. The overall distribution of geographical origins is, as expected, quite skewed: The six most frequent countries of citizenship among all speakers -- United States (24.6\%), France (13.2\%), Germany  (9.5\%), Russia  (8.5\%), United Kingdom (7.4\%) and Italy (7.2\%) -- comprise more than 70\% of all. The evaluation for women yields a picture similar to the overall trend: almost the same countries appear under the top six, comprising more than 72\% of all talks by women: United States (28\%) and France (18.3\%) are the most frequent countries, followed by Germany (8.3\%) and United Kingdom (6.9\%). Russia and Italy are less strong than in the distribution of all ICM speakers, with 5.5\% and 4.6\%, respectively. Italy, in fact, is not among the top six countries for women speakers, ranking at position seven, slightly behind Israel with a share of 5\% among female speakers. Israel clearly marks the most notable difference in the comparison of origins, with an overall share of less than 2\% compared to 5\% among women speakers. 

Almost all congresses show a certain regional focus, manifested in the composition of the organizing committees as well as the origin of the invited speakers (see Tab. \ref{tab:1} in Appendix). For instance, more than 32\% of all speakers at the ICM-1900 in Paris were French; 45\% of all speakers at the congress in Heidelberg in 1904 were German; United Kingdom was the country of citizenship for more than 26\% of all speakers in 1912 in Cambridge (UK). The focus on 'local' speakers, evident in subsequent congresses as well, is especially pronounced for congresses that took place in countries with overall strong representation such as Russia at ICM-1966 in Moscow or the United States at ICM-1986 in Berkeley. Nevertheless, many other countries managed to foster the representation of their scientists. Brasil, for instance, accounts for 6\% of all speakers at the latest congress in Rio de Janeiro. 
A look at the origins of women shows a rather inconclusive picture. A similar trend is expressed in a more drastic way: when a congress took place in a country with an overall strong representation, often all invited women originated from the host country. At most of these congresses, however, this corresponded to only one or two invited women. Nevertheless, at some congresses with a stronger representation of women such as the ICM in 1928 in Bologna or in 1986 in Berkeley, the host country reached a share of 40 to 50\% among the invited women. On the other hand, there are 12 ICMs to which women were invited but none originated from the host country.

\subsection{Topics not balanced}
The individual character of the congresses is reflected by the diversity of names chosen for the themed sections, summing up to more than 180 different titles. While some section names such as \textit{Numerical Methods, Numerical Mathematics} or \textit{Numerical Methods and Computing} obviously belong together, this is much less the case with topics such as \textit{History of Mathematics, Logic and Foundations} and \textit{Mathematics Education}: at various ICMs, two or even all three of them had been combined together into one section, making it impossible to study these topics individually without intensive manual work. 
The division into few, rather broad sections was typical for the early ICMs. The ICM-1928 in Bologna, for instance, combined in the first section talks on \textit{Algebra, Arithmetics and Analysis}, and had only one additional section on pure mathematics, mainly for talks in \textit{Geometry}. \textit{Elementary Mathematics}, \textit{Didactical Questions}, and \textit{Mathematical Logic} were grouped into one section, and \textit{Philosophy and History of Mathematics} into another one. On the contrary, recent ICMs feature more than 20 themed sections, providing a better granularity to analyze the share of women speakers by their  field of research. 

The data by the IMU \cite{ICM} shows certain inconsistencies, in particular for some early ICMs. For instance, Section 1 at ICM-1928 was named \textit{Analysis} instead of \textit{Arithmetics, Algebra, Analysis}. Furthermore, 15 talks from 1966 were assigned to the section 1/2 hr. report. The distinction between 1 hour and 1/2 hour is basically the distinction between plenary talks and section talks, and while most of the 1/2 hour talks were  assigned to a themed section, these 15 remained for unknown reasons. 

The partially unclean data for ICMs before 1970, and the fact that except for the congresses in 1928 and 1932 the presence of women before the 1980s is highly scattered, have prompted us to restrict our description of women's participation by sections to the period 1970--2018. 

\begin{figure}[b]
\centering
\includegraphics[width=\textwidth]{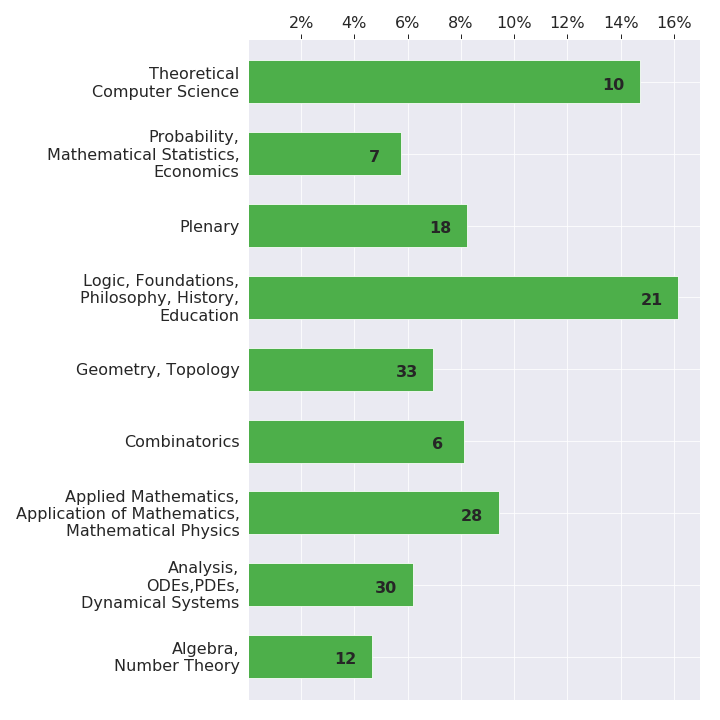}
\caption{Percentage of talks by women per group of sections since 1970. The numbers inside the bars correspond to the total number of talks by women.}
\label{fig:sections}      
\end{figure}

We have grouped the section names manually into the following eight groups (number of talks per section since ICM-1970 in parentheses): \textit{Logic, Foundations, Philosophy, History and Education} (155); \textit{Applied Mathematics, Applications of Mathematics, Mathematical Physics} (327); \textit{Probability, Mathematical Statistics, Economics} (129); \textit{Analysis, ODEs, PDEs, Dynamical Systems} (519); \textit{Algebra and Number Theory} (271); \textit{Theoretical Computer Science} (78); \textit{Geometry and Topology} (515); \textit{Combinatorics} (80). In addition, there were a total of 237 \textit{Plenary} talks since 1970. There are two further categories of talks: the \textit{ICM Abel Lecture} given to the winner of the Abel Prize (2 talks, both by men) and the \textit{ICM Emmy Noether Lecture}\footnote{The Emmy Noether Lecture honors women who have made fundamental and sustained contributions to the mathematical sciences. Since 2006, this lecture is a permanent ICM feature, since 2014, a special commemorative plaquette is conferred to every ICM Emmy Noether Lecturer.
} (6 talks by women). We have omitted both of them from the plot in Fig. \ref{fig:sections}. 

As shown in Fig. \ref{fig:sections}, sections concerning \textit{Algebra and Number Theory} have the least 
proportion of talks by women (<5\%), closely followed by sections dealing with \textit{Probability, Mathematical Statistics and Economics} (<6\%). The two largest groups of sections, \textit{Analysis, ODEs, PDEs, Dynamical Systems}, and \textit{Geometry and Topology}, which together comprise more than 1,000 talks since 1970, have both less than 7\% of women speakers on average. On the positive end we see two rather small section groups, \textit{Logic, Foundations, Philosophy, History, Education} (>16\%) and \textit{Theoretical Computer Science} (>14\%). The \textit{Plenary} section, which is supposed to contain the most prominent congress talks, contains 8.2\% talks by women and is hence slightly above the average of 7.3\% since 1970. 

Further investigation needs to be carried out to understand the unequal distribution among topics. In any case, there is no conclusive correlation between this distribution and the representation of women authors in the respective fields. The distribution of authorships and authors across classes of the Mathematics Subject Classification (MSC) 2010 has been carried out in \cite{MST} based on data from zbmath.org, one of the two main services for bibliographic information in Mathematics. Fig. 10 in \cite{MST} shows the amount of women authors as a heatmap across MSC classes. It shows, for instance, that an above-average number of women publish in the field of Statistics and their share in Probability theory is close to the overall average. At the same time, the group of ICM sections related to \textit{Statistics, Probability and Economics} is, with less than 6\%, almost at the very bottom of the scale. Likewise, women authors are very well represented in most MSC classes related to \textit{Analysis, PDEs, ODEs and Dynamical Systems}, in particular in relatively large fields like PDEs and ODEs. However, the respective group of ICM sections, while being the largest in terms of the total number of talks, has very few  talks by women. 

In \cite{Sad}, Cora Sadosky notes that the distribution of women across non-plenary sections suggests an ``invisible quota system'', leading to at most one woman per congress and section. She continues saying that ``it seemed as if the selection panels, although aware enough to consider women candidates, felt that they had fulfilled their duty when the first one accepted.'' The first exceptions to this pattern occur in the section \textit{Mathematical Aspects of Computer Science} featuring two talks by women in 1990 (Lenore Blum and Shafi Goldwasser) and in 1998 (Joan Feigenbaum and Toniann Pitassi). Also, the section related to \textit{Teaching and Popularization of Mathematics} often contained more than one talk by a woman. In fact, most topics featured two talks by a woman at some congresses. In some rather exceptional cases, even more than two women were invited, such as in 2014 in \textit{Combinatorics} or \textit{Mathematics in Science and Technology}, or in 2018 in the \textit{Geometry} section. Sections with an approximately equal distribution of women and men are extremely rare; the small section \textit{Mathematics Education and Popularization of Mathematics} at the last ICM in Rio de Janeiro with two talks by women and one talk by a man seems to be the only example so far. Nevertheless, the presence of at most a single woman among the section speakers as criticized by Sadosky remains by far the usual practice.

\section{Discussion}
In our study we have used the complete list of plenary and section speakers invited to ICMs from its beginning in 1897 through the most recent congress in 2018 in Rio de Janeiro. We have combined this list with demographic information from Wikidata and section names of the talks. We have further enriched our data, in particular by inferring the gender. Our data is provided at \cite{Mih} and can be used for further research.

The participation of women in the International Congresses of Mathematicians has increased over the course of the past 121 years. However, as shown in Fig. \ref{fig:longitudinal_dev}, their share does not show a clear trend over time. Instead, the inclusion of women among ICM speakers shows a clear peak in the late 1920s and early 1930s of the last century, followed by a long period of almost complete absence, until the start of a continuous positive development in the late 1980s that persists until today. As noted by Izabella Laba with respect to the equitable representation of women speakers at ICM-2014, ``Compared to the proportion of women among tenured and tenure-track faculty at research universities, the group from which ICM invited speakers usually hail and therefore a more appropriate benchmark, it does not look far off'' \cite{Lab}.

We have analyzed the relation of the speakers' gender with their countries of citizenships and the topics of the sections in which they presented their talks. As expected, most speakers came from a rather small set of countries, while many countries and even entire continents were barely represented at all. The distribution of citizenships with respect to gender, however, does not show significant difference. More than 70\% of all women and men speakers are or were citizens of countries whose territory today corresponds to the United States, France, Russia, Germany, United Kingdom and Italy. 
We have further shown that the selection of speakers reflects a certain regional focus, yielding, in most cases, a noticeably higher share of speakers originating from the congress’s host country. This effect is particularly pronounced for countries that are overall well represented among ICM speakers. For women, however, this tendency is not so clear and would require further investigation. Cost and ease of travel as well as the language are likely to have played a role in the decisions to invite a speaker, in particular in the early congresses. Nevertheless, it would be interesting to explore the possible relations between the nationalities of the members of the Programme Committee and of the sectional committees with the nationalities of the invited women (and men).

The individual congresses show high variance in the arrangements of talks into sections, making a grouping of sections inevitable. We have arranged the non-plenary talks since ICM-1970 into eight large groups and studied the distribution of women across these. We have found that the share of women ranges between less than 5\% in sections related to Algebra and Number Theory to more than 16\% in sections on Logic, Foundations, Philosophy, History and Education. Talks by women sum up to 6\% of all talks in the section group Analysis, ODEs, PDEs, Dynamical Systems which is by far the largest one. These effects are not consistent with the representation of women as authors in the respective fields. However, a correlation on this level might not be very likely anyway since ICM invitations are very rare and, in some sense, singular (for women). An analysis of gender and topics among authors in highly prestigious journals or among tenured and tenure-track faculty at research universities would hence presumably constitute a more suitable set in order to understand whether the fluctuations between ICM sections might partially be explained by a lack of `suitable' women in the respective fields. It should, however, be noted that the breakdown by sections leads to a rather small number of individuals per group, which is more prone to variation. 

It would be interesting to consider the longitudinal development of women’s participation in other conferences in mathematics (and other fields) as well. Smaller conferences might constitute a much bigger issue in terms of inclusion since, as claimed in \cite{Lab}, there are no large committees overseeing them. In recent years, (women) scientists in STEM fields have proposed the formulation of a Bechdel Test, a measure of women’s representation in fiction (movies, comics, video games etc.), for scientific workshops. A movie would pass the Bechdel test if (1) it features at least two named women, (2) who talk to each other, (3) about something besides a man. An analogous test for scientific conferences could require (1) at least two female invited speakers, (2) who are asked questions by female audience members, (3) about their research\footnote{http://openscience.org/a-bechdel-test-for-scientific-workshops/}. As noted by various scientists in social media channels, this form of test is rarely passed by conferences in STEM fields; even among the recent ICMs, many sections would not pass this test either. Such a test, while measuring only a basic level of inclusion of women and despite being far from creating a `critical mass' in the respective conferences, would yield a first understanding of women’s participation in mathematical conferences across fields and over time.

Manifold factors have played a role in the longitudinal evolution of the number of ICM lecturers at ICM. Of particular importance for the sustained positive evolution in the last decades was and is the establishment of various associations of women in mathematics and their efforts to increase the visibility of their female colleagues in the field. As suggested by Elena Resmerita and Carola-Bibiane Schönlieb,  the current convenors of the European Women in Mathematics (EWM), it is crucial to keep highlighting contributions of women mathematicians by continuing to showcase contributions of women in mathematics, ``raise the profile of women mathematicians, volunteer to serve on committees of international mathematical associations and mathematical award committees, nominate female colleagues or encourage others to nominate them, and overall help to build a scientific atmosphere without boundaries.''\footnote{Statement by the EWM convenors Elena Resmerita and Carola-Bibiane Schönlieb in reaction to the absence of women among the Fields medalists in 2018. https://www.europeanwomeninmaths.org/ewm-statement-about-icm-2018}

\section*{Acknowledgements}
This research was conducted as part of the project `A Global Approach to the Gender Gap in Mathematical, Computing, and Natural Sciences: How to Measure It, How to Reduce It?', supported by the International Science Council and several scientific unions. We thank June Barrow-Green, Fulvia Furinghetti, Catherine Goldstein and  Annette Vogt for providing important insight and expertise from History of Mathematics. We further thank Lucia Santamaria for suggestions that greatly improved the manuscript. We gratefully acknowledge the support of the International Mathematical Union for granting us access to their data on ICM Plenary and Invited Speakers.

\section*{Appendix}
\addcontentsline{toc}{section}{Appendix}
%
%

%
\begin{table}
\caption{Percentage of all speakers and all women speakers, respectively, whose country of citizenship equals the host country of the respective congress.}
\label{tab:1} 
\begin{tabular}{p{2cm}p{4cm}p{3cm}p{3cm}}
\hline\noalign{\smallskip}
year & host country & speakers from host country (\%)  & women speakers from host country (\%)  \\
\noalign{\smallskip}\svhline\noalign{\smallskip}
1897 & Switzerland & 3.1 & 0 \\
1900 & France & 32.4 & 0 \\
1904 & Germany & 45.2 & 0 \\
1908 & Italy & 32.7 & 100 \\
1912 & United Kingdom & 26.5 & 100 \\
1920 & France & 31.4 & 0 \\
1924 & Canada & 2.3 & 0 \\
1928 & Italy & 25 & 50 \\
1932 & Switzerland & 7.4 & 8.3 \\
1936 & Norway & 6.1 & 0 \\
1950 & United States of America & 36.6 & 0 \\
1954 & Netherlands & 2.5 & 0 \\
1958 & United Kingdom & 6.9 & 0 \\
1962 & Sweden & 2.9 & 0 \\
1966 & Russia & 30.4 & 100 \\
1970 & France & 8.1 & 0 \\
1974 & Canada & 2.3 & 0 \\
1978 & Finland & 0.5 & 0 \\
1983 & Poland & 3 & 0 \\
1986 & United States of America & 39.0 & 40 \\
1990 & Japan & 10.1 & 0 \\
1994 & Switzerland & 1.7 & 0 \\
1998 & Germany & 7.6 & 0 \\
2002 & China & 3.3 & 5.3 \\
2006 & Spain & 3.8 & 0 \\
2010 & India & 7.7 & 3.6 \\
2014 & South Korea & 0 & 0 \\
2018 & Brazil & 5.9 & 12.5 \\
\noalign{\smallskip}\hline\noalign{\smallskip}
\end{tabular}
\vspace*{-12pt}
\end{table}

%
%
%

\end{document}